\documentclass[12pt, reqno]{amsart}
\usepackage{amssymb,a4}
\usepackage{amsmath,amsthm,amsfonts,amssymb}
\def\cen{\centerline}

\def\z{\mathbb{Z}}

\def\dim{\hbox{dim}}
\def\ad{\hbox{ad}}

\def\a{\alpha}

\def\sg{\sigma}

\def\gl{\frak{gl}}
\newfont{\df}{eufm10}

\def\sg{\sigma}

\def\der{\hbox{Der}\,}

\def\fsl{\frak{sl}\,}

\def\stl{\frak{stl}\,}

\def\vp{\varphi}

\def\ot{\otimes}

\def\de{\delta}
\def\De{\Delta}
\def\dim{\hbox{\rm dim}\,}

\def\str{\hbox{\rm str}}
\def\hom{\hbox{\rm Hom}\,}
\def\en{\hbox{\rm End}\,}
\def\ker{\hbox{Ker}\,}

\def\ad{\hbox{\rm ad}\,}

\def\ot{\otimes}

\def\mg{{\bf \frak g}}
\def\dmg{\dot{\mg}}

\hoffset
\voffset
\title{Leibniz superalgebras graded by finite root systems}

\author{Naihong Hu}
\address{Department of Mathematics, East China Normal
University, Shanghai, 200062, P.R. China}
\email{nhhu@gmath.ecnu.edu.cn}

\author{Dong Liu}
\address{Department of Mathematics, Huzhou Teachers College, Zhejiang Huzhou, 313000, China}
\email{liudong@hutc.zj.cn, corresponding arthor}

\author{Linsheng Zhu}
\address{Department of Mathematics, Changshu Institute of
Technology, Jiangsu Changshu, 215500, China}
\email{lszhu@cslg.edu.cn}

\thanks{}

\begin{document}
\maketitle

\begin{abstract}
The structure of Lie algebras, Lie superalgebras and Leibniz
algebras graded by finite root systems has been studied by several
researchers since 1992. In this paper, we study the structure of
Leibniz superalgebras graded by finite root systems, which gives
an approach to study various classes of Leibniz superalgebras.
\smallskip

\end{abstract}

\newtheorem{theo}{Theorem}[section]
\newtheorem{defi}[theo]{Definition}
\newtheorem{lemm}[theo]{Lemma}
\newtheorem{coro}[theo]{Corollary}
\newtheorem{prop}[theo]{Proposition}

\section{Introduction}

In \cite{Lo1}, J.-L. Loday introduces a non-antisymmetric version of
Lie algebras, whose bracket satisfies the Leibniz relation and
therefore is called {\it Leibniz algebra}. The Leibniz relation,
combined with antisymmetry, is a variation of the Jacobi identity,
hence Lie algebras are anti-symmetric Leibniz algebras. In
\cite{Lo2}, Loday also introduces an {\sl associative} version of
Leibniz algebras, called {\it associative dialgebras}, equipped with
two binary operations, $\vdash$ and $\dashv$, which satisfy the five
relations (see the axiom (Ass) in Section 2). These identities are
all variations of the associative law, so associative algebras are
dialgebras for which the two products coincide. The peculiar point
is that the bracket $[a, b]=:a\vdash b-b\dashv a$ defines a (left)
Leibniz algebra which is not antisymmetric, unless the left and
right products coincide. Hence dialgebras yield a commutative
diagram of categories and functors
\begin{eqnarray*}
{\bf Dias}&\stackrel{-}{\to}& {\bf Leib}\\
\downarrow&&\downarrow\\
{\bf Ass}&\stackrel{-}{\to}& {\bf Lie}
\end{eqnarray*}

Recently super associative diagebras and Leibniz superalgebras were
studied in \cite{D}, \cite{LH3} and \cite{AAO}, etc.. The structure
of Lie algebras, Lie superalgebras and Leibniz algebras graded by
finite root systems was studied in several papers (\cite{BM},
\cite{BZ}, \cite{BE1}---\cite{BE3}, \cite{BEM}, \cite{ABG1},
\cite{ABG2}, \cite{LH2}, etc.). In this paper we determine the
structure of Leibniz superalgebras graded by the root systems of the
basic classical simple Lie superalgebras.

The paper is organized as follows. In Section 2, we recall some
notions of Leibniz superalgebras and superdialgebras. In Section 3,
we give the definition and  some properties of Leibniz superalgebras
graded by finite root systems. In Section 4 and Section 5, we mainly
determine the structure of Leibniz superalgebras graded by the root
systems of types $A(m, n)$ and $C(n)$, $D(m, n)$, $D(2, 1; \a)$,
$F(4)$, $G(2)$.

Throughout this paper, $K$ denotes a field of characteristic 0. $A$
denotes a unital associative super associative diagebra over $K$.

\section{ASSOCIATIVE SUPER DIALGEBRAS AND LEIBNIZ SUPERALGEBRAS}

We recall some notions of super associative diagebras and Leibniz
superalgebras and their (co)homology as defined in \cite{D} and
\cite{LH3}.

\subsection{Associative super associative diagebras.}

\begin{defi} \cite{Lo2} An {\it associative dialgebra} $D$ over $K$ is a $K$-vector space
with two operations $\dashv, \vdash:D\ot D\to D$, called left and
right products, satisfying the following five axioms:
$$\begin{cases} a\dashv(b\dashv c)=(a\dashv b)\dashv c=a\dashv(b\vdash c),\\
        (a\vdash b)\dashv c=a\vdash(b\dashv c),\\
         (a\vdash b)\vdash c=a\vdash (b\vdash c)=(a\dashv b)\vdash c.  \end{cases} \leqno(Ass)$$
\end{defi}

Obviously an associative dialgebra is an associative algebra if
and only if $a\dashv b=a\vdash b=ab$.

An {\it associative super associative diagebra} over $K$ is a
$\z_2$-graded $K$-vector space $D$ with two operations $\dashv,
\vdash:D\ot D\to D$, satisfying the axiom (Ass) and
$$D_{\sg}\dashv D_{\sg'}, D_{\sg}\vdash D_{\sg'} \subset
D_{\sg+\sg'},\quad \forall\sg, \sg'\in\z_2.$$

A super associative diagebra is called commutative if $a\vdash
b=b\dashv a$ for all $a, b\in D$.  A super associative diagebra is
called unital if it is given a specified bar-unit: an element $1\in
D$ which is a unit for the left and right products only on the
bar-side, that is $1\vdash a=a=a\dashv 1$, for any $a\in D$. Denote
by {\bf SDias, SAss} the categories of super associative diagebras
or associative superalgebras over $K$ respectively. Then the
category {\bf SAss} is a full subcategory of {\bf SDias}.

\noindent{\bf Examples}. 1. Obviously, a super associative diagebra
is an associative superalgebra if $a\dashv b=a\vdash b=ab$. An
associative dialgebra is a trivial super associative diagebra.

2. {\it Super differential dialgebra.} Let $(A=A_{\bar0}\oplus
A_{\bar1}, d)$ be a differential associative super
algebra($|d|=\bar0$). So by hypothesis, $d(ab)=(da)b+adb$ and
$d^2=0$. Define left and right products on $A$ by the formulas
$x\dashv y=xdy$ and $x\vdash y=(dx)y$. Then $A$ equipped with these
two products is a super associative diagebra.

3. {\it Tensor product.} If $D$ and $D'$ are two super associative
diagebras, then the tensor product $D\ot D'$ is also a super
dialgebra with
$$(a\ot a')\star(b\ot b')=(-1)^{|a'||b|}(a\star b)\ot (a'\star b')\eqno(2.1)$$ for $\star=\dashv, \vdash.$

For instance, $M_{m+n}(D):=M_{m+n}(K)\ot D$ is a super associative
diagebra if $D$ is a super associative diagebra and $M_{m+n}(K)$ is
the superalgebra of all $(m+n)\times (m+n)$-matrices over $K$.

4. The {\it free associative super associative diagebra} (see
\cite{LH3} in details) on an $\z_2$-graded vector space $V$ is the
dialgebra $Dias(V)=T(V)\ot V\ot T(V)$ equipped with the induced
$\z_2$-gradation.

\subsection{Leibniz superalgebra}

\begin{defi} \cite{D}  A Leibniz superalgebra is a $\z_2$-graded vector space
$L=L_{\bar 0}\oplus L_{\bar 1}$ over a field $K$ equipped with a
$K$-bilinear map $[-,-]: L\times L\to L$ satisfying

$$[L_{\sg}, L_{\sg'}]\subset L_{\sg+\sg'},\quad \forall \sg, \sg'\in\z_2$$
and the Leibniz superidentity
$$[[a, b], c]= [a, [b, c]]-(-1)^{|a||b|}[b, [a, c]], \quad \forall \;a, \,b, \,c\in L.\eqno(2.2)$$

\end{defi}

Obviously, $L_{\bar0}$ is a Leibniz algebra. Moreover any Lie
superalgebra is a Leibniz superalgebra and any Leibniz algebra is
a trivial Leibniz superalgebra. A Leibniz superalgebra is a Lie
superalgebra if and only if
$$[\,a, b\,]+(-1)^{|a||b| }[b, a]=0, \quad \forall\; a, b\in L.\eqno(2.3)$$

\noindent{\bf Examples.} 1. Let ${\mg}$ be a Lie superalgebra, $D$
be a unital commutative dialgebra, then $\mg\ot D$ with Leibniz
bracket $[x\ot a, y\ot b]=[x, y]\ot(a\vdash b)$ is a Leibniz
superalgebra. Let ${\mg}$ be a  basic classical simple Lie
superalgebra which is not of type $A(n, n), \, n\ge 1$, then
$\tilde\mg=\mg\ot D\oplus \Omega_D^1$ with the bracket
$$[x\ot a, y\ot b]=[x, y]\ot(a\vdash b) +(x, y)b\dashv da, \quad \forall a, b\in D, x, y\in \mg,\eqno(2.4)$$
$$[\Omega_D^1,\  \tilde\mg]=0\eqno(2.5)$$
is also a Leibniz superalgebra, where $(-, -)$  is an even
invariant bilinear form of $\mg$, $\Omega_D^1$ is defined in
\cite{L} (also see \cite{LH3}). In fact it is the universal
central extension of $\mg\ot D$ (see \cite{LH3} in details).

2. Tensor product. Let $\mg$ be a Lie superalgebra,  then the
bracket
$$[x\ot y, a\ot b]=[[x, y], a]\ot b+(-1)^{|a||b|}a\ot [[x, y], b]\eqno(2.6)$$ defines a Leibniz superalgebra structure on the vector space $\mg\ot\mg$
(see \cite{KP} for that in Leibniz algebras case).

3. The {\it general linear Leibniz superalgebra} $\gl(m, n, D)$ is
generated by all $n\times n$ matrices with coefficients from a
dialgebra $D$, and $m, n\ge 0, n+m\ge2$ with the bracket
$$[E_{ij}(a), E_{kl}(b)]=\de_{jk}E_{il}(a\vdash b)-(-1)^{\tau_{ij}\tau_{kl}}\de_{il}E_{kj}(b\dashv a),\eqno(2.7)$$ for all $a, b\in D$.

Clearly, $\gl(m, n, D)$ is a Leibniz superalgebra. If $D$ is an
associative superalgebra, then  $\gl(m, n, D)$ becomes a Lie
superalgebra.

By definition, the {\it special linear Leibniz superalgebra} with
coefficients in $D$ is
$$\fsl(m, n, D):=[\,\gl(m, n, D), \gl(m, n, D)\,].$$ Notice that if $n\ne m$ the Leibniz superalgebra $\fsl(m, n, D)$ is simple.

The special linear Leibniz superalgebra $\fsl(m, n, D)$ has
generators $E_{ij}(a), 1\le i\ne j\le m+n, a\in D$, which satisfy
the following relations:
\begin{eqnarray*}
&&[E_{ij}(a), E_{kl}(b)]=0\  \hbox{if}\ i\ne l\ \hbox{and}\  j\ne k;\\
&&[E_{ij}(a), E_{kl}(b)]=E_{il}(a\vdash b)\  \hbox{if}\ i\ne l\ \hbox{and}\  j= k;\\
&&[E_{ij}(a), E_{kl}(b)]=-(-1)^{\tau_{ij}\tau_{kl}}E_{kj}(b\dashv
a)\ \hbox{if}\ i= l\ \hbox{and}\  j\ne k,
\end{eqnarray*}

4.  The {\it Steinberg Leibniz superalgebra} $\stl(m, n, D)$
(\cite{L}) is a Leibniz superalgebra generated by symbols
$u_{ij}(a)$, $1\le i\ne j\le n$, $a\in D$, subject to the
relations
\begin{eqnarray*}
&&v_{ij}(k_1a+k_2b)=k_1v_{ij}(a)+k_2v_{ij}(b),\  \hbox{ for } \ a, b\in D, \ k_1, k_2\in K;\\
&&[v_{ij}(a), v_{kl}(b)]=0,\  \hbox{ if }\ i\ne l\ \hbox{ and } \  j\ne k;\\
&&[v_{ij}(a), v_{kl}(b)]=v_{il}(a\vdash b)\  \hbox{ if } \ i\ne l\ \hbox{and } \  j= k;\\
&&[v_{ij}(a), v_{kl}(b)]=-(-1)^{\tau_{ij}\tau_{kl}}v_{kj}(b\dashv
a)\  \hbox{ if } \ i= l\ \hbox{and} \  j\ne k,
\end{eqnarray*}
where $1\le i\ne j\le m+n$, $a\in D$. It is clear that the last
two relations make sense only if $m+n\ge3$. See \cite{L} and
\cite{LH4} for more details about the Steinberg Leibniz
superalgebra.

We also denote by {\bf SLeib, SLie} the categories of Leibniz
superalgebras and Lie superalgebras over $K$ respectively.

For any associative super associative diagebra $D$, if we define
$$[x, y]=x\vdash y-(-1)^{|x||y|}y\dashv x,\eqno(2.8)$$
then $D$ equipped with this bracket is a Leibniz superalgebra. We
denoted it by $ D_L$. The canonical morphism $D\to D_{L}$ induces
a functor $(-):$ {\bf SDias}$\to${\bf SLeib}.

\noindent{\bf Remark.} For a super associative diagebra $D$, if we
define
$$[x, y]=x\dashv y-(-1)^{|x||y|}y\vdash x,\eqno(2.9)$$ then $(D, [
, ])$ is a right Leibniz superalgebra in the sense of \cite{LP}.

For a Leibniz superalgebra $L$, let $L_{LS}$ be the quotient of
$L$ by the ideal generated by elements $[x, y]+(-1)^{|x||y|}[y,
x]$, for all $x, y\in L$. It is clear that $ L_{LS}$ is a Lie
superalgebra. The canonical epimorphism $: L\to L_{LS}$ is
universal among the maps from $L$ to Lie superalgebras. In other
words the functor $(-)_{LS}:$ {\bf SLeib}$\to${\bf SLie} is left
adjoint to $inc:$ {\bf SLie}$\to${\bf SLeib}.

Moreover we have the following commutative diagram of categories
and functors
\begin{eqnarray*}
{\bf SDias}&\stackrel{-}{\to}& {\bf SLeib}\\
\downarrow&&\downarrow\\
{\bf SAss}&\stackrel{-}{\to}& {\bf SLie}
\end{eqnarray*}

As in the Leibniz algebra case, the {\it universal enveloping super
associative diagebra} (\cite{LH3}) of a Leibniz superalgebra $L$ is
$$Ud(L):=(T(L)\ot L\ot T(L))/\{[x, y]-x\vdash y+(-1)^{|x||y|}y\dashv x|x, y\in L\}.$$

\begin{prop} \cite{LH3}
 The functor $Ud:{\bf SLeib}\to {\bf SDias}$ is left adjoint to the functor $-:{\bf SDias}\to {\bf SLeib}$. \hfill $\rule[-.23ex]{1.0ex}{2.0ex}$
\end{prop}

Let $L$ be a Leibniz superalgebra. We call a $\z_2$-graded space $M=M_{\bar0}\oplus M_{\bar1}$ a left $L$-module if there is a bilinear map:
$$[-,-]: L\times M\to M $$
satisfying the following three axioms
$$[[x, y], m]=[x, [y, m]]-(-1)^{|x||y|}[y, [x, m]],$$
for any $m\in M$ and $x,\, y\in L$. In this case we also say $\vp:
L\to \en_KM,\ \vp(x)(m)=[x, m]$ is a left representation of $L$.

Clearly, if $L$ is a Lie superalgebra, then $M$ is just the left $L$-module in the Lie superalgebra case.

Suppose that $L$ is a Leibniz superalgebra over $K$. For any $z\in
L$, we define $\ad z\in \hbox{End}_kL$ by
$$\ad z(x)=[z, x], \quad\forall x\in L.\eqno(2.10)$$
It follows (2.2) that
$$\ad z([x, y])=[\ad z(x), y]+(-1)^{|z||x|}[x, \ad z(y)]\eqno(2.11)$$
for all $x, y\in L$. This says that $\ad z$ is a super derivation
of degree $|z|$ of $L$. We also call it the inner derivation of
$L$.

Similarly we also have the definition of general super derivation
of a Leibniz superalgebra $L$. By definition a super derivation of
degree $s, s\in \z_2$ of $L$ is an endomorphism
$\mu\in\hbox{End}_sL$ with the property
$$\mu([a, b])=[\mu(a), b]+(-1)^{s|a|}[a, \mu(b)].$$
We denote by Inn($L$), Der$(L)$ the sets of all inner derivations,
super derivations of L respectively. They are also Leibniz
superalgebras.

For a Lie superalgebra $L$, $HL^1(L, M)=H^1(L, M)=\der(L,
M)/\hbox{Inn}(L, M)$, where $H^1(L, M)$ (resp. $HL^1(L, M)$) denotes
the first cohomology groups of $L$ in the category of Lie (resp.
Leibniz) superalgebras(see \cite{Lo1}).

\subsection{ Leibniz algebras graded by finite root
systems}

Now we recall some notions of Leibniz algebras graded by finite
root systems defined as in \cite{LH2}.
\begin{defi} \cite{LH2}
A Leibniz algebra $L$ over a field $K$ of characteristic 0 is
graded by the $($reduced$)$ root system $\Delta$ or is
$\Delta$-graded if

$(1)$ L contains as a subalgebra a finite-dimensional simple Lie
algebra $\dmg=H\oplus\bigoplus_{\a\in\Delta}\dmg_{\a}$ whose root
system is $\Delta$ relative to a split Cartan subalgebra $H=\dmg_0$;

$(2)$ $L=\bigoplus_{\a\in\Delta\cup\{0\}}L_{\a}$, where
$L_{\a}=\{x\in L\mid \ad h(x)=-[x, h]=\a(h)x, \forall h\in H\}$
for $\a\in \Delta\cup\{0\}$; and

$(3)$ $L_0=\sum_{\a\in\Delta}[L_\a, L_{-\a}]$.
\end{defi}

\noindent{\bf Remarks.}

1. The conditions for being a $\De$-graded Leibniz algebra imply
that $L$ is a direct sum of finite-dimensional irreducible Leibniz
representations of $\dmg$ whose highest weights are roots, hence
are either the highest long root or high short root or are $0$.

2. If $L$ is $\Delta$-graded then $L$ is perfect (i.e. $[L, L]=L$).
Indeed the result follows from $L_\a=[L_\a, H]$ for all
$\a\in\Delta$ and (3) above.

In \cite{LH2}, the structure of Leibniz algebras graded by the root
systems of types $A, D, E$ was determined by using the methods in
\cite{BM}. In fact we have

\begin{theo} \cite{LH2}
 Let $L$ be a Leibniz algebra over $K$ graded by the root system $\Delta$ of type $X_l\,(l\ge2)\, (X_l=A_l, D_l, E_l)$ .

$(1)$ If $X_l=A_l, l\ge 3$, then there exists a unital associative
$K$-dialgebra $R$ such that $L$ is centrally isogenous with $\fsl(n,
R)$, .

$(2)$  If $X_l=A_l, l=2$, then there exists a unital alternative
$K$-dialgebra $R$ such that $L$ is centrally isogenous with the
Steinberg Leibniz algebra $\stl(n, R)$, where $\stl(n, R)$ defined
in \cite{L} (also see the example 4 in Section 2.2).

$(3)$  If $X_l=D_l\, (l\ge 4), E_l\, (l=6, 7, 8)$, then there exists
a unital associative commutative $K$-dialgebra $A$ such that $L$ is
centrally isogenous with the Leibniz algebra $\dmg\ot R$.

\end{theo}

\noindent{\bf Remark.} Two perfect Lie algebras $L_1$ and $L_2$ are
called {\it centrally isogenous} if they have the same universal
central extension (up to isomorphism).

\section{Leibniz superalgebras graded by finite root systems}

In this section, we shall study the structure of Leibniz
superalgebras graded by finite root systems.

Motivated the definitions of Lie superalgebras and Leibniz
algebras graded by finite root systems defined in
\cite{BE1}-\cite{BE3} and \cite{LH2}, we give the following
definition.

\begin{defi}
A Leibniz superalgebra $L$ over a field $K$ of characteristic 0 is graded by the $($reduced$)$ root system $\Delta$ or is $\Delta$-graded if

$(1)$ L contains as a subsuperalgebra a finite-dimensional split simple basic Lie superalgebra $\mg=H\oplus\bigoplus_{\a\in\Delta}\mg_{\a}$ whose root system $\Delta$ is relative to a split Cartan subalgebra $H=\mg_0$;

$(2)$ $L=\bigoplus_{\a\in\Delta\cup\{0\}}L_{\a}$, where $L_{\a}=\{x\in L\mid \ad h(x)=[h, x]=\a(h)x, \forall h\in H\}$ for $\a\in \Delta\cup\{0\}$; and

$(3)$ $L_0=\sum_{\a\in\Delta}[L_\a, L_{-\a}]$.
\end{defi}

\noindent{\bf Remarks.}

1. If $L$ is $\Delta$-graded then $L$ is perfect. Indeed  the
result follows from $L_\a=[H, L_\a]$ for all $\a\in\Delta$ and (3)
above.

2. The Steinberg Leibniz superalgebra $\stl(m, n, D)\,(m\ne n)$ is graded by the root system of type $A_{m-1, n-1}$.

We would like to view $L$ as a left $\mg$-module in order to
determine the structure of $L$. The following result plays a key
role in examining $\De$-graded Leibniz superalgebras. It follows
from the Lemma 2.2 in \cite{BE3}.
\begin{lemm}
Let $L$ be a $\De$-graded Leibniz superalgebra, and let $\mg$ be
its associated split simple basic Lie superalgebra. Then $L$ is
locally finite as a module for $\mg$. \hfill $\quad
\rule[-.23ex]{1.0ex}{2.0ex}$
\end{lemm}

As a consequence, each element of a $\De$-graded Leibniz
superalgebra $L$, in particular each weight vector of $L$ relative
to the Cartan subalgebra $H$ of $\mg$, generates a
finite-dimensional $\mg$-module. Such a finite-dimensional module
has a $\mg$-composition series whose irreducible factors have
weight which are roots of $\mg$ or $0$. Next we determine which
finite-dimensional irreducible $\mg$-modules have nonzero weights
that are roots of $\mg$.

Throughout this paper we will identify the split simple Lie
superalgebras $\mg$ of type $A(m, n),\, m>n\ge0$, with the special
linear Lie superalgebra $\fsl(m+1, n+1)$. For simplicity of
notation, set $p=m+1, q=n+1$, so $\mg=\fsl(p, q), \, p>q\ge1$.

\begin{prop}\cite{BE1}, \cite{BE2}
Let $\mg$ be a split simple Lie superalgebra of type $A(m, n)$, with $m\ge n\ge 0, m+n\ge1$, or $C(n)$, $D(m, n)$, $D(2, 1; \a)\ (\a\not\in\{0, -1\})$, $F(4)$, $G(2)$. The only finite-dimensional irreducible left $\mg$-modules whose weights relative the Cartan subalgebra of diagonal matrices (modulo the center if necessary) are either roots or 0 are exactly the adjoint and the trivial modules ( possibly with the parity changed).
\end{prop}

\begin{prop}\cite{BE1}, \cite{BE2}
Let $\mg$ be a split simple Lie superalgebra of type $A(m, n)\ (m>n\ge 0)$ or $C(n)$, $D(m, n)$, $D(2, 1; \a)\ (\a\not\in\{0, -1\})$, $F(4)$, $G(2)$, with split Cartan subalgebra $H$. Assume $V$ is a locally finite left $\mg$-module satisfying

(i) $H$ acts semisimply on $V$;

(ii) any composition factor of any finite-dimensional right submodule of $V$ is isomorphic to the adjoint representation $\mg$ or to a trivial representation.

Then $V$ is completely reducible. \hfill$\quad
\rule[-.23ex]{1.0ex}{2.0ex}$

\end{prop}

\begin{prop}\cite{BE1}
Let $\mg$ be a split simple Lie superalgebra of type $A(m, n)$, with $m>n\ge 0$. Then $\hom_\mg(\mg\ot\mg, \mg)$ is two-dimensional and spanned by the Lie (super) bracket and by the map given by $(x, y)\mapsto x\ast y=xy+(-1)^{|x||y|}yx-{2\over m-n}\str(xy)I$, for any $x, y\in \mg=\fsl(m+1, n+1)$. $\quad \rule[-.23ex]{1.0ex}{2.0ex}$
\end{prop}

%\newpage

\section{The structure of the $A(m, n)$-graded Leibniz superalgebras $(m>n)$}

The results of Section 2 show that any $A(m, n)$-graded Leibniz
superalgebra $L$ with $m>n\ge0$ is the direct sum of adjoint and
trivial modules (possibly with a change of parity) for the grading
subalgebra $\mg$. After collecting isomorphic summands, we may
suppose that there are superspaces $A=A_{\bar0}\oplus A_{\bar1}$ and
$D=D_{\bar0}\oplus D_{\bar 1}$ so that $L=(\mg\ot A)\oplus D$, and a
distinguished element $1\in A_{\bar 0}$ which allows us to identify
the grading subalgebras $\mg$ with $\mg\ot1$. Observe first that $D$
is a subsuperalgebra of $L$ since it is the (super) centralizer of
$\mg$.

To determine the multiplication on $L$, we may apply the same type
of arguments as in \cite{BE2}. Indeed, fix homogeneous basis
elements $\{a_i\}_{i\in I}$ of $A$ and choose $a_i, a_j, a_k$ with
$i, j, k\in I$, we see that the projection of the product $[\mg\ot
a_i, \mg\ot a_j]$ onto $\mg\ot a_k$ determines an elements of
$\hom_\mg(\mg\ot\mg, \mg)$, which is spanned by the Leibniz
supercommutator. Then there exist scalars $\xi_{i, j}^k$ and
$\theta_{i, j}^k$ so that
$$[x\ot a_i, y\ot a_j]|_{\mg\ot A}=(-1)^{|a_i||y|}\left([x, y]\ot (\sum_{k\in I}\xi_{i, j}^ka_k)+x\ast y\ot (\sum_{k\in I}\theta_{i, j}^ka_k)\right)$$

Define $\circ: A\times A\to A$ by $a_i\circ a_j=2\sum_{k\in
I}\xi_{i, j}^ka_k$, and $[ , ]: A\times A\to A$ by $[a_i,
a_j]=2\sum_{k\in I}\theta_{i, j}^ka_k$ and extending them
bilinearly, we obtain two products ``$\circ$'' and ``[ , ]'' on
$A$.

Taking into account that $\hom_\mg(\mg\ot \mg, K)$ is spanned by
the supertrace, we see that there exists a bilinear form $\langle,
\rangle:A\times A\to D$ and an even bilinear map $D\times A\to
A:(d, a)\to da$ with $d1=0$ such that the multiplication in $L$ is
given by

$$[f\ot a, g\ot b]=(-1)^{|a||g|}\left([f, g]\ot {1\over 2}a\circ b+f\ast g\ot {1\over2}[a, b]+\str(fg)\langle a, b\rangle\right),\eqno(4.1)$$
$$[d, f\ot a]=(-1)^{|d||f|}f\ot da,\eqno(4.2)$$
for homogeneous elements $f, g\in\mg,\  a, b\in A, d\in D$.
Additionally, $1\circ a=2a$ and $[1, a]=0$ for all $a\in A$.

There are two unital multiplications $\dashv$ and $\vdash$ on $A$
such that
$$a\circ b=a\vdash b+(-1)^{|a||b|}b\dashv a,\eqno(4.3)$$
$$[a, b]=a\vdash b-(-1)^{|a||b|}b\dashv a,\eqno(4.4)$$
for any homogeneous elements $a, b\in A$. Moreover by setting $a=1$, we have $1\vdash b=b\dashv 1=b$,
so the dialgebra $A$ is unital.

Now the Jacobi superidentity $[[z_1, z_2], z_3]=[z_1, [z_2,
z_3]]-(-1)^{|z_1||z_2|}[z_2, [z_1, z_3]]$, when specialized with
homogeneous elements $d_1, d_2\in D$ and $f\ot a\in \mg\ot A$,
shows that  $\phi: D\to \en_K(A): \phi(d)(a)=da$, is a left
representation of the Leibniz superalgebra $D$. When it is
specialized with homogeneous elements $d\in D$ and $f\ot a, g\ot
b\in \mg\ot A$, we obtain
$$[d, [f\ot a, g\ot b]]=[[d, f\ot a], g\ot b]+(-1)^{|d|(|f|+|a|)}[f\ot a, [d, g\ot b]],$$
and using (4.1) and (4.2), we see that this is same as:
\begin{eqnarray*}
&&(-1)^{|d|(|f|+|g|)}(-1)^{|a||g|}\left([f, g]\ot {1\over 2}d(a\circ b)+f\ast g\ot {1\over 2}d([a, b])+\str(fg)[d, \langle a, b\rangle]\right)\\
&=&(-1)^{|d||f|}(-1)^{|d|+|a||g|}\left([f, g]\ot {1\over 2}da\circ b)+f\ast g\ot {1\over 2}[da, b]+\str(fg)\langle da, b\rangle\right)\\
&+&(-1)^{|d|(|f|+|g|+|a|)}(-1)^{|a||g|}\left([f, g]\ot {1\over 2}a\circ db+f\ast g\ot {1\over 2}[a, db]+\str(fg)\langle a, db\rangle\right).\\
\end{eqnarray*}
When $f=E_{1, 2}$ and $g=E_{2, 1}$, the elements $[f, g]$ and $f\ast g$ are linearly independent and $\str(fg)=1$. Hence we have

(i) $d(a\circ b)=(da)\circ b+(-1)^{|d||a|}a\circ (db),$

(ii) $d([a, b])=[da, b]+(-1)^{|d||a|}[a, db],$

(iii) $[d, \langle a, b\rangle]=\langle da,
b\rangle+(-1)^{|d||a|}\langle a, db\rangle$,

for any homogeneous $d\in D$ and $a, b\in A$.

Items (i) and (ii) can be combined to give that
$$\phi\ \hbox{is a left representation as superderivations}:\phi: D\to \der_K(A).\eqno(4.5)$$
While (iii) says that
$$\langle\  ,\, \rangle\ \hbox{is invariant under the action of } D.\eqno(4.6)$$

For $f\ot a, g\ot b, h\ot c\in \mg\ot A$, the Jacobi superidentity is equivalent to the following two relations ($\mg\ot A$ and $D$ components):
\begin{eqnarray*}
(\star)&&\Big(\str([f, g]h)\langle {1\over 2}a\circ b, c\rangle+\str((f\ast g)h)\langle {1\over 2}[a, b], c\rangle\Big)\\
&-&\Big(\str(f[g, h])\langle a, {1\over 2}b\circ c\rangle+\str(f(g\ast h))\langle a, {1\over 2}[b, c]\rangle\Big)\\
&+&(-1)^{(|f||g|+|a||b|)}\Big(\str(g[f,  h])\langle b, {1\over 2}a\circ c\rangle+\str(g(f\ast h))\langle b, {1\over 2}[a, c]\rangle\Big)\\
&=&0.
\end{eqnarray*}
and
\begin{eqnarray*}
(\star\star)&&\Big([[f, g], h]\ot {1\over 4}(a\circ b)\circ c+[f, g]\ast h\ot {1\over 4}[a\circ b, c]\\
&&+[f\ast g,  h]\ot{1\over 4}[a, b]\circ c+(f\ast g)\ast h\ot {1\over 4}[[a, b], c]+\str(fg)h\ot \langle a, b\rangle c\Big)\\
&-&\Big([f, [g, h]]\ot {1\over 4}a\circ (b\circ c)+f\ast [g, h]\ot {1\over 4}[a, b\circ c]\\
&&+ [f,  g\ast h]\ot {1\over 4}a\circ [b, c]+f\ast (g\ast h)\ot {1\over 4}[a, [b, c]]+(-1)^{|a|(|g|+|h|)}\str(gh)[f\ot a, \langle b, c\rangle]\Big)\\
&+&(-1)^{|f||g|+|a||b|}\Big([g, [f, h]]\ot {1\over 4}b\circ (a\circ c)+g\ast [f, h]\ot {1\over 4}[b, a\circ c]\\
&&+[g,  f\ast h]\ot {1\over 4}b\circ [a, c]+g\ast (f\ast h)\ot {1\over 4}[b, [a, c]]\\
&&+(-1)^{|b|(|f|+|h|)} \str(fh)[g\ot b,\langle a, c\rangle]\Big)\\
&=&0.
\end{eqnarray*}

The formula $(\star)$ can be written as
\begin{eqnarray*}
&&\str(fgh)\big(\langle a\vdash b, c\rangle-\langle a, b\vdash c\rangle-(-1)^{|a|(|b|+|c|)}\langle b, c\dashv a\rangle\big)\\
&-&(-1)^{|g||h|+|a||b|}\str(fhg)\big(\langle b\dashv a, c\rangle-\langle b, a\vdash c\rangle-(-1)^{|b|(|a|+|c|)}\langle a, c\dashv b\rangle\big)\\
&=&0.
\end{eqnarray*}

Then we have
$$\langle a\vdash b, c\rangle=\langle a, b\vdash c\rangle+(-1)^{|a|(|b|+|c|)}\langle b, c\dashv a\rangle,\eqno(4.7)$$
$$\langle a\vdash b, c\rangle=\langle a\dashv b, c\rangle.\eqno(4.8)$$

The formula $(\star\star)$ can be written as:
\begin{eqnarray*}
(\star\star\star)&&fgh\ot\big((a\vdash b)\vdash c-a\vdash (b\vdash c)\big)\\
&-&(-1)^{|f||g|+|a||b|}gfh\ot\big((b\dashv a)\vdash c-b\vdash (a\vdash c)\big)\\
&+&(-1)^{(|f|+|g|)|h|+(|a|+|b|)|c|}hfg\ot\big((c\dashv a)\dashv b-c\dashv (a\vdash b)\big)\\
&+&(-1)^{(|h|+|g|)|f|+(|c|+|b|)|a|}ghf\ot\big((b\vdash c)\dashv a-b\vdash (c\dashv a)\big)\\
&+&(-1)^{|h||g|+|c||b|}fhg\ot\big((a\vdash c)\vdash b-a\vdash (c\dashv b)\big)\\
&-&(-1)^{(|h||g|+|f||g|+|f||h|+|c||b|+|a||b|+|a||c|}hgf\ot\big((c\dashv b)\dashv a-c\dashv (b\dashv a)\big)\\
&-&\str(fg)\left({1\over m-n}h\ot [[a, b], c]-h\ot \langle a, b\rangle c\right)\\
&+&\str(gh)\left({1\over m-n}f\ot[a, [b, c]]-(-1)^{|a|(|g|+|h|)}[f\ot a, \langle b, c\rangle]\right)\\
&-&(-1)^{|f||g|+|a||b|}\str(fh)\left({1\over m-n}g\ot [b, [a, c]]-(-1)^{|b|(|f|+|h|)}[g\ot b,\langle a, c\rangle]\right)\\
&-&{1\over m-n}\str(fgh)I\ot\Big([a\vdash b, c]-[a, b\vdash c]-(-1)^{|a|(|b|+|c|)}[b, c\dashv a]\Big)\\
&+&(-1)^{|f||g|+|a||b|}{1\over m-n}\str(gfh)I\ot\Big([b\dashv a, c]-(-1)^{|b|(|a|+|c|)}[a, c\dashv b]-[b, a\vdash c]\Big)\\
&=&0.
\end{eqnarray*}

Set $f=E_{12}, g=E_{23}, h=E_{31}$ in $(\star\star\star)$, if $m\ge 2$
then by $|f|=|g|=|h|=\bar 0$ and the independent of $E_{11}, E_{22}, E_{33}, I$, we have:
$$(a\vdash b)\vdash c=a\vdash(b\vdash c),\eqno(4.9)$$
$$(b\vdash c)\dashv a=b\vdash (c\dashv a),\eqno(4.10)$$
$$(c\dashv a)\dashv b=c\dashv(a\vdash b).\eqno(4.11)$$

Similarly by setting $f=E_{31}, g=E_{23}, h=E_{12}$, we can obtain
$$(c\dashv b)\dashv a=c\dashv(b\dashv a),\eqno(4.12)$$
$$(b\dashv a)\vdash c=b\vdash (a\vdash c),\eqno(4.13)$$
so $A$ is an associative super associative diagebra.

If $m=1$, then $|f|=\bar 0$ and $|g|=|h|=\bar 0$. The expression in $(\star\star\star)$ is a linear combination of $E_{11}, E_{22}, E_{33}$ with coefficents in $A$. By direct calcualtion we also obtain that $A$ is associative.

Then $(\star\star\star)$ becomes
\begin{eqnarray*}
&-&\str(fg)\left({1\over m-n}h\ot [[a, b], c]-h\ot \langle a, b\rangle c\right)\\
&+&\str(gh)\left({1\over m-n}f\ot[a, [b, c]]-(-1)^{|a|(|g|+|h|)}[f\ot a, \langle b, c\rangle]\right)\\
&-&(-1)^{|f||g|+|a||b|}\str(fh)\left({1\over m-n}g\ot [b, [a, c]]-(-1)^{|b|(|f|+|h|)}[g\ot b,\langle a, c\rangle]\right)=0
\end{eqnarray*}

Then $$\langle a, b\rangle c={1\over m-n}[[a, b], c]\eqno(4.14)$$
and
$$[f\ot a, \langle b, c\rangle]={1\over m-n}f\ot [a, [b, c]]\eqno(4.15)$$ since $|g|+|h|=\bar 0$ if $\str(gh)\ne 0$.

In this way, we have arrived at our main Theorem. The last sentence
in it is a consequence of condition (3) in Definition 3.1.

\begin{theo}
Assume $L=\mg\ot A\oplus D$ is a superalgebra over a field $K$ of
characteristic $0$ where $\mg=\fsl(m+1, n+1), m>n\ge 0$, $A$ is a
unital super associative diagebra, and $D$ is a Leibniz
superalgebra, and with multiplication as in (4.1) and (4.2). Then
$L$ is a Leibniz superalgebra if and only if

(1) $A$ is an associative super associative diagebra,

(2) $D$ is a Leibniz subsuperalgebra of $L$ and $\phi: D\to \der_K(A)(\phi(d)a=da)$ is a left representation of $D$ as superderivations on the dialgebra $A$,

(3) $[d, \langle a, b\rangle]=\langle da, b\rangle+(-1)^{|d||a|}\langle a, db\rangle$,

(4) (4.7), (4.8) and (4.14), (4.15) hold,

for any homogeneous elements $d\in D$ and $a, b, c\in A$.

Moreover, the $A(m, n)$-graded Leibniz superalgebras (for $m>n\ge 0$) are exactly these superalgebras with the added constraint that
$$D=\langle A, A\rangle.$$
\end{theo}

\noindent{\bf Remark.} Let $A$ be any unital associative super
dialgebra. Then $\ad_{[A, A]}$ is a subsuperalgebra of $\der_K(A)$
(it is a Lie superalgebra). Consider the Leibniz superalgebra
$$\frak L(A):=(\mg\ot A)\oplus \ad_{[A, A]},$$
with $\mg=\fsl(m+1, n+1)(m>n\ge0)$, with multiplication given by
(4.1) and (4.2) in place of $D$ and with $\langle a,
b\rangle={1\over m-n}\ad_{[A, A]}$ for any $a, b\in A$. Then Theorem
4.1 shows that $\frak L(A)$ is an $A(m, n)$-graded Leibniz
superalgebra. Moreover for any $A(m, n)$-graded Leibniz superalgebra
$L$ with coordinate super associative diagebra $A$, Theorem 4.1
implies that $L/Z(L)\cong \frak L(A)$. Thus $L$ is a cover of $\frak
L(A)$.

\begin{coro}
The $A(m, n)$-graded Leibniz superalgebras (for $m>n\ge 0$) are
precisely the Leibniz superalgebras which are centrally isogeneous
to the Leibniz superalgebra $\fsl(m+1, n+1, A)$ for a unital
associative super associative diagebra $A$.
\end{coro}

\noindent{\bf Remarks.}  1.  The situation when $m=n$ in $A(m, n)$
is much more involved than the case of $m\ne n$, due to the fact
that the complete reducibility in Proposition 3.4 no longer valid
in this case. However, using the similar consideration as above
and that in \cite{BEM}, we can also obtain the following result.

\begin{theo}
Let L be a Leibniz superalgebra graded by $\mg$, which is a split
simple classical Lie superalgebra of type $A(n, n) (n>1)$. Then
there exists a unital associative super associative diagebra $D$
such that it is centrally isogenous to
$$\frak{sl}(n+1, D)=[\frak{gl}(n+1)\ot D, \frak{gl}(n+1)\ot D].$$
\end{theo}

2. For a unital associative dialgebra $A$, the universal central
extension of the Leibniz superalgebra $\fsl(m, n, A)$ with $m\ne n$
and $m+n\ge 3$ has been shown to be the Steinberg Leibniz supergebra
$\stl(m, n, A)$ in \cite{L}.

\section{The structure of $\De$-graded Leibniz superalgebras of other types}

It follows from Proposition 3.4 that every Leibniz superalgebra
graded by the root system $C(n)$, $D(m, n)$, $D(2, 1; \a)\
(\a\not\in\{0, -1\})$ , $F(4)$, or $G(2)$ decomposes as a
$\mg$-module into a direct sum of adjoint modules and trivial
modules. The next general result describes the structure of Leibniz
superalgebras $L$ having such decompositions, which is essentially
same as the Lemma 4.1 in \cite{BE2}.

\begin{lemm}
Let $L$ be a Leibniz superalgebra over $K$ with a subsuperalgebra $\mg$, and assume that under the adjoint action of $\mg$, $L$ is a direct sum of

$(1)$ copies of the adjoint module $\mg$,

$(2)$ copies of the trivial module $K$.

Assume that

$(3)$ $\dim\hom_\mg(\mg\ot\mg, \mg)=1$ so that $\hom_\mg(\mg\ot\mg, \mg)$ is spanned by $x\ot y\to [x, y]$.

$(4)$ $\hom_\mg(\mg\ot\mg, K)=K\kappa$, where $\kappa$ is even,
non-degenerate and supersymmetric, and the following conditions
hold:

$(5)$ There exist $f, g\in \mg_{\bar0}$ such that $[f, g]\ne0$ and $\kappa(f, g)\ne 0$,

$(6)$ There exist $f, g, h\in \mg_{\bar0}$ such that $[f, h]=[g, h]=0$ and $\kappa(f, h)=\kappa(g, h)=0\ne \kappa(f, g)$,

$(7)$ There exist $f, g, h\in \mg_{\bar0}$ such that $[[f, g], h]=0\ne [[g, h], f]$.

Then there exist superspaces $A$ and $D$ such that $L\cong(\mg\ot
A)\oplus D$ and

$(a)$ $A$ is a unital supercommutative associative super associative
diagebra;

$(b)$ $D$ is a trivial $\mg$-module and is a Leibniz superalgebra;

$(c)$ Multiplication in $L$ is given by

$$[f\ot a, g\ot b]=(-1)^{|a||g|}\Big([f, g]\ot (a\vdash b)+\kappa(f, g)\langle a, b\rangle\Big),\eqno(5.1)$$
$$[d, f\ot a]=(-1)^{|d||f|}f\ot da,\eqno(5.2)$$
for homogeneous elements $f, g\in\mg, a, b\in A, d\in D$, where

(i) $D$ is a Leibniz subsuperalgebra of $L$ and $\phi: D\to \der_K(A)(\phi(d)a=da)$ is a left representation of $D$ as superderivations on the dialgebra $A$ with $\langle A, A\rangle\subset \ker\phi$,

(ii) $[d, \langle a, b\rangle]=\langle da, b\rangle+(-1)^{|d||a|}\langle a, db\rangle$, in particular, $\langle A, A\rangle$ is an ideal of $D$,

(iii) $\langle a\vdash b, c\rangle=\langle a, b\vdash c\rangle+(-1)^{|a|(|b|+|c|)}\langle b, c\dashv a\rangle$ and $\langle a\vdash b, c\rangle=\langle a\dashv b, c\rangle.$

Conversely, the conditions above are sufficient to guarantee that
a superspace $L=\mg\ot A\oplus D$ satisfied (a)--(c) is a Leibniz
superalgebra.
\end{lemm}
\noindent{\bf Proof.} When a Leibniz superalgebra $L$ is a direct
sum of copies of adjoint modules and trivial module for $\mg$,
then after collecting isomorphic summands, we may assume that
there are superspaces $A=A_{\bar 0}\oplus A_{\bar 1}$ and
$D=D_{\bar 0}\oplus D_{\bar 1}$ so that $L=\mg\ot A\oplus D$.
Suppose that such a superalgebra $L$ satisfies conditions
(1)---(4). Note that $D$ is a super subalgebra of $L$.

Fixing homogeneous basis $\{a_i\}_{i\in I}$ of $A$ and choose
$a_i, a_j, a_k$ with $i, j, k\in I$, we see that the projection of
the product $[\mg\ot a_i, \mg\ot a_j]$ onto $\mg\ot a_k$
determines an elements of $\hom_\mg(\mg\ot\mg, \mg)$, which is
spanned by the Leibniz supercommutator. Then there exist scalars
$\xi_{i, j}^k$ and $\theta_{i, j}^k$ so that
$$[x\ot a_i, y\ot a_j]|_{\mg\ot A}=(-1)^{|a_i||y|}[x, y]\ot (\sum_{k\in I}\xi_{i, j}^ka_k).$$

Defining $\dashv, \vdash: A\times A\to A$ by $a_i\vdash
a_j=\sum_{k\in I}\xi_{i, j}^ka_k=(-1)^{|a_i||a_j|a_j\vdash a_i}$
and extending it bilinearly, we have a supercommutative dialgebra
structure on $A$.

By similar arguments, there exist bilinear maps $A\times A\to D,
(a, b)\mapsto \langle a, b\rangle\in D$, and  $D\times A\to A, (d,
a)\mapsto da\in A$, such that the multiplication in $L$ is as in
(c).

First the Jacobi superidentity, when specialized with homogeneous
elements $d_1, d_2\in D$ and $f\ot a\in \mg\ot A$, and then with
$d\in D$ and $f\ot a, g\ot b\in \mg\ot A$ will show that
$\phi(d)a=da$ is a representation of $D$ as superderivations of
$A$. We assume next that $f, g$ are taken to satisfy (5). Then for
homogeneous elements $d\in D,\, a, b\in A$, the Jacobi
superindentity gives the condition $[d, \langle a,
b\rangle]=\langle da, b\rangle+(-1)^{|d||a|}\langle a, db\rangle$.

The Jacobi superindentity, when specialized with homogeneous
elements $f\ot a, g\ot b, h\ot c\in\mg\ot A$ and $f, g, h$ as
condition (7), gives that $\langle A, A\rangle$ is contained in
the kernel of $\phi$.

Finally the Jacobi superindentity, when specialized with homogeneous
elements $f\ot a, g\ot b, h\ot c\in\mg\ot A$ and $f, g, h$ as
condition (8), gives the (iii) and associativity of dialgebra $A$
(see Section 4 or \cite{BE1}).

The converse is a simple computation. \hfill
$\rule[-.23ex]{1.0ex}{2.0ex}$

From \cite{BE1}, we see that $\mg$ satisfies the above condition if
$\mg$ is a split simple classical Lie superalgebra of type $C(n)
(\ge3)$, $D(m, n)(m\ge2, n\ge1)$, $D(2, 1; \a), \a\in K, \a\ne 0,
1$, $F(4)$ or $G(2)$. Then we have the following structure theorem
for Leibniz superalgebras graded by the root system of type $\mg$
(compare Theorem 5.2 in \cite{BE1}).

\begin{theo}
Let $\mg$ be a split simple classical Lie superalgebra of type
$C(n)\, (\ge3)$, $D(m, n)\,(m\ge2, n\ge1)$, $D(2, 1; \a), \,\a\in
K,\, \a\ne 0, 1$, $F(4)$ or $G(2)$, then by \cite{BE1} we have
$\dim\hom_\mg(\mg\ot\mg, \mg)=\dim\hom_\mg(\mg\ot\mg, K)=1$. Let $L$
be a $\De$-graded Leibniz superalgebra of type $\mg$, then there
exists a unital supercommutative associative super associative
diagebra $A$ and a $K$-superspace $D$ such that $L\cong(\mg\ot
A)\oplus D$. Multiplication in $L$ is given by

$$[f\ot a, g\ot b]=(-1)^{|a||g|}\Big([f, g]\ot (a\vdash b)+\kappa(f, g)\langle a, b\rangle\Big),\eqno(5.3)$$
$$[d, L]=0,\eqno(5.4)$$
for homogeneous elements $f, g\in\mg, \, a, b\in A, \, d\in D$,
where $\kappa( , )$ is a fix even nondegenerate supersymmetric
invariant bilinear form on $\mg$ and $\langle , \rangle: A\ot A\to
D$ is a $K$-bilinear and satisfies the following conditions:

$(1)$ $[d, \langle a, b\rangle]=\langle da, b\rangle+(-1)^{|d||a|}\langle a, db\rangle$,

$(2)$ $\langle a\vdash b, c\rangle=\langle a, b\vdash
c\rangle+(-1)^{|a|(|b|+|c|)}\langle b, c\dashv a\rangle$,\quad
$\langle a\vdash b, c\rangle=\langle a\dashv b, c\rangle$,

and $D$ is a Leibniz subsuperalgebra of $L$ and $\phi: D\to \der_K(A)(\phi(d)a=da)$ is a left representation of $D$ as superderivations on the dialgebra $A$.
\end{theo}
\noindent{\bf Proof.} The only point left is the proof of the
centrality of $D$. Condition (3) of Definition 3.1 forces $D=\langle
A, A\rangle$, which by Lemma 5.1 is contained in $\ker\phi$.
Therefore $D=\langle A, A\rangle$ is abelian and centralizes $\mg\ot
A$, hence it is central.
 \hfill $\rule[-.23ex]{1.0ex}{2.0ex}$

As the similar arguments in \cite{BE1}, we also have

\begin{coro}
Let $\mg$ be a split simple classical Lie superalgebra of type $C(n)
(\ge3)$, $D(m, n)\,(m\ge2, n\ge1)$, $D(2, 1; \a),\, \a\in K, \,\a\ne
0, 1$, $F(4)$ or $G(2)$ and $L$ be a $\De$-graded Leibniz
superalgebra of type $\mg$, then there exists a unital
supercommutative associative super associative diagebra $A$ such
that $L$ is a covering of the Leibniz superalgebra $\mg\ot A$.
\hfill $\rule[-.23ex]{1.0ex}{2.0ex}$

\end{coro}

Let $\mg$ be a split simple classical Lie superalgebra of type $C(n)
(\ge3)$, $D(m, n)\,(m\ge2,\, n\ge1)$, $D(2, 1; \a),\, \a\in K,\,
\a\ne 0, 1$, $F(4)$ or $G(2)$, $A$ a unital associative commutative
dialgebra, the universal central extension of the Leibniz
superalgebra $\mg\ot A$ has been studied in \cite{LH3}.

\noindent{\bf Remark.} Such researches for the cases of $B(m, n)$ and $A(1, 1)$ are
more complicated (see \cite{BE3}, \cite{BEM}), which depends on the results of (non)associative diagebras.

\vskip45pt \centerline{\bf ACKNOWLEDGMENTS}

\vskip15pt

N. Hu
is supported in part by by the NNSF (Grant 10431040, 10271047), the
TRAPOYT and the FUDP from the MOE of China, the SRSTP from the
STCSM, the Shanghai Priority Academic Discipline from the SMEC. 
D. Liu is supported by the NNSF (No. 11071068, 10701019), the ZJNSF(No. D7080080,Y6100148),
Qianjiang Excellence Project (No. 2007R10031), the "New Century 151
Talent Project" (2008) and the "Innovation Team Foundation of the
Department of Education"  (No. T200924) of Zhejiang Province.

\vskip30pt
\def\refname{\cen{\normalsize\bf REFERENCES}}

\end{document}